\documentclass[12pt]{article}
\usepackage{amsmath,amsfonts,amsthm}

\def\hybrid{\topmargin 0pt      \oddsidemargin 0pt
        \headheight 0pt \headsep 0pt
        \voffset=-0.5cm
        \textwidth 6.5in        
        \textheight 9in         
        \marginparwidth 0.0in
        \parskip 5pt plus 1pt   \jot = 1.5ex}

\hybrid

\def\beq{\begin{equation}}
\def\eeq{\end{equation}}
\def\bea{\begin{eqnarray}}
\def\eea{\end{eqnarray}}
\def\p{\partial}
\def\G{\Gamma}
\def\g{\gamma}
\def\s{\sigma}

\def\C{{\cal C}}

\def\a{\alpha}

\def\A{{\cal A}}

\def\D{{\cal D}}

\def\h{{\cal H}}
\def\L{{\cal L}}

\def\M{{\cal M}}
\def\MS{{\cal M}_{g,\,k}^{<n>}}
\def\PS{{\cal P}_{g,\,k}^{<n>}}

\def\Poo{{\cal P}_{g,1}^{<1>}}
\def\MO{{\cal M}_{g,\,k}^{(\,n\,)}}
\def\Mo{{\cal M}_{g,1}^{(\,n\,)}}
\def\Moo{{\cal M}_{g,1}^{(\,1\,)}}
\def\Mon{{\cal M}_{g,1}^{<n>}}
\def\MOO{{\cal M}_{g,\,k}^{(\,m\,)}}
\def\Mom{{\cal M}_{g,1}^{(\,m\,)}}

\def\W{{\cal W}}
\def\SP{{\cal S}}

\def\wt{\widetilde}

\def \matrix #1 {\left(\begin{array}{cc} #1 \end{array}\right)}

\newtheorem{theo}{Theorem}[section]

\newtheorem{cor}[theo]{Corollary}
\newtheorem{lem}[theo]{Lemma}

\theoremstyle{definition}
\newtheorem{dfn}[theo]{Definition}
\newtheorem{rem}[theo]{Remark}

\begin{document}

\title{Real normalized differentials and Arbarello's conjecture}
\author{I.Krichever
\thanks{Columbia University, New York, USA, Landau
Institute for Theoretical Physics and Kharkevich
Institute for Problems of Information Transmission, Moscow, Russia ; e-mail:
krichev@math.columbia.edu. Research is supported in part
by The Ministry of Education and Science of the Russian Federation (contract 02.740.11.5194).}}

\date{}

\maketitle

\begin{abstract}
\noindent
Using meromorphic differentials with real periods, we prove Arbarello's conjecture:
any compact complex cycle of dimension $g-n$ in the moduli space $\M_g$ of smooth genus $g$ algebraic curves must intersect the locus of curves having a Weierstrass point of order at most $n$.

\end{abstract}

\section{Introduction}

In \cite{arb} Arbarello constructed a natural filtration
\beq\label{arbstrat}
\W_2\subset\W_3\subset\cdots\subset\W_{g-1}\subset\W_g=\M_g
\eeq
of the moduli space $\M_g$ of smooth Riemann surfaces of genus $g$, whose stratum $\W_n$ is the locus of smooth genus $g$ algebraic curves having a Weierstrass point of order at most $n$
(i.e. the locus of curves on which there exists a meromorphic function with one pole of order at most $n$). He conjectured that: {\it any compact complex cycle in $\M_g$ of dimension $g-n$ must intersect the stratum $\W_n$}. Since $\W_2$ is the locus of hyperelliptic curves, which is affine, Arbarello's conjecture implies that {\it $\M_g$ does not contain complete (complex) subvarieties of dimension greater than $g-2$}. That was proved later by Diaz in \cite{diaz} with the help of a variant of Arbarello's stratification. Another modification of Arbarello's stratification was used by Looijenga who proved that {\it the tautological classes of degree greater than $g-2$ vanish in the Chow ring of $\M_g$}. The later implies Diaz' result (the Hodge class $\lambda_1$ is ample on $\M_g$, and thus for any complete $d$-dimensional subvariety $X\subset\M_g$ we would have $\lambda_1^d\cdot X>0$, while $\lambda_1^{g-1}=0$, as a tautological class).

The main goal of this paper is to prove Arbarello's conjecture, which until now has remained open.
(A highly non-trivial nature of this problem has found its "explanation" in recent work \cite{arbmond}, where it was shown that open parts of strata $\W_n$ are almost never affine). Our proof uses certain constructions of the Whitham perturbation theory of integrable systems \cite{kr-av,kr-tau}, further developed and clarified in \cite{kp1,kp2}. These constructions has already found their applications in topological quantum field theories (WDVV equations) and $N=2$ supersymmetric gauge theories \cite{gmmmk} (see also \cite{brad} and references therein). Applications of these constructions to a study of geometry of the moduli spaces of curves was initiated in author's joint works with S. Grushevsky (\cite{kr-grush1,kr-grush2}).

In \cite{kr-grush1} we gave a new proof of Diaz' theorem, and in \cite{kr-grush2} we proved vanishing of certain tautological classes. The both results had been known, but their new proofs were the evidence that further developments of the Whitham constructions, which was the main goal of \cite{kr-grush1,kr-grush2}, can make them instrumental for a variety of algebraic-geometric problems.

A notion of {\it real-normalized} meromorphic differentials is central in the Whitham theory. By definition a real normalized meromorphic differential is a differential whose periods over any cycle on the curve are real. The power of this notion is that on any algebraic curve and for each fixed set of "singular parts" at the marked labeled points there exists a unique real normalized differential
having prescribed singularities.

In \cite{kr-grush1} the real normalized differentials of the third kind, i.e. differentials with two simple poles were used. The statement of Arbarello's conjecture is about curves with one puncture. That's why it was tempting for the author to use real normalized differentials of the second kind for its proof. It turned out, that a naive attempt to use arguments analogous to ones in \cite{kr-grush1} runs almost immediately into a obstacle due to a non-compactness of the space of the singular parts defining such differentials uniquely. This obstacle is parallel to ones in the original attempts to prove Arbarello's conjecture (see details in the book \cite{har-mor}).

In the next section necessary extensions and partial compactification of the previously
known constructions are presented. Namely, we define a foliation structure on the space of real normalized differentials of the second kind with poles of orders {\it at most} $n_\a+1$ at marked points.
In early works the foliation structure was defined on the moduli space of real normalized differentials with {\it fixed} order of poles. It turned out that, as before, each leaf
of the foliation is (locally) smooth complex subvariety. At the end of Section 2 we show
that the foliation structure on the moduli space of real normalized differential
descents to a foliation structure on factor space of these differentials  considered modulo action
of the multiplicative group of {\it positive real numbers}.

In Section 3 we introduce  additional tool needed for the proof of the main theorem --- that is a notion of cycles "dual" to zeros of real normalized differentials. We prove that "dual cycles"
span the homology group $H_1(\G,\mathbb Z)$. It seems for the author that the construction of this cycles is of independent interest and we return to it elsewhere. The proof of Arbarello's conjecture is given in Section 4.

{\sl Acknowledgments.} The author would like to thank Sam Grushevsky for the numerous valuable comments, clarifications, and suggestions.

\section{Foliations defined by real normalized differentials}

Let $\MO,\ n=(n_1,\ldots,n_k)$ be the moduli space of smooth genus $g$ algebraic curves $\G$ with fixed {\it singular part of a pole of order $n_\a+1$ with no residue} in the neighborhood of labeled marked distinct points  $p_\a\in \G,\, \a=1,\ldots,k$. We recall that explicitly choosing a singular part of a pole of order at $n_\a+1$ near $p_\a$ means on a small neighborhood of $p_\a$ choosing: (i) a coordinate $z_\a$ such that $z_\a(p_\a)=0$; (ii) polynomial $R_\a$ of the form $R_\a=\sum_{i=0}^{n_\a}r_{\a,\,i}z_\a^{-i-1}$, and identifying pairs $(z_\a,R_\a)$ and $(w_\a, R'_\a)$ if $R'\,dw_\a=R_\a\,dz_\a+O(1)\,dz_\a,\, \, w_\a=w_\a(z_\a)$. The coefficient $r_{\a,\,0}$ is the residue of the singular part, i.e. for singular parts with no residues $r_{\a,\,0}=0$.

Non-degeneracy of the imaginary part of Riemann matrix of $b$-periods of normalized holomorphic differential on a smooth genus $g$ algebraic curve $\G$ implies that:

\begin{lem}\label{1} For any fixed singular parts of poles with pure imaginary residues, there exists a unique meromorphic differential $\Psi$, having prescribed singular part at $P_\a$ and such that all its periods on $\G$ are real, i.e.
\beq\label{realnorm}
{\rm Im}\, \left(\oint_c \Psi\right)=0, \ \ \forall\  c\in H^1(\G,\mathbb Z)
\eeq
\end{lem}
\noindent
(for detailed proof see Proposition 3.4 in \cite{kr-grush1}).

\begin{rem} Although throughout the paper we consider only  real normalized differentials of the second kind (i.e. having no residues at poles), most of the constructions can be easy extended to the case of meromorphic differential with residues, which have to be pure imaginary. For the first time, real normalization as defining property of quasi-momentum differentials in the spectral theory of linear operators with quasi-periodic coefficients was introduced in \cite{kr-real} and \cite{kr-av} (where they were called absolute normalized).
\end{rem}
Let $\MS$ be the moduli space of smooth genus $g$ algebraic curves $\G$ with fixed {\it non-trivial} real normalized meromorphic differential $\Psi$ of the second kind, having poles of order at most $n_\a+1$ at the labeled points $p_\a\in \G$. Lemma \ref{1} identifies $\MS$ with the moduli space of curves with fixed non trivial set of singular parts of poles of orders {\it at most} $n_{\a}+1$ with no residues at the marked points. The later is the total space of rank $|n|:=\sum_\a{n_a}$ vector bundle over $\M_{g,k}$ with zero section removed. Therefore,
the identification of a real normalized differential with its singular part defines on $\MS$ a complex structure.

\medskip
Our next goal is to introduce on $\MS$ the structure of a foliation.
Periods of the differential $\Psi$ define a cohomology class $\Pi\in H^1(\G, \mathbb R)$. Let $\nabla$ be the Gauss-Manin connection on the Hodge bundle over $\M_{g,k}$, whose fiber over $(\G,p_\a)$ is $H^1(\G,\mathbb Z)$.
Then the equation
\beq
\nabla_X \, \Pi=0
\eeq
considered as an equation for the tangent vector $X\in T(\MS)$ defines at every point of $\MS$ a subspace of the tangent space. A distribution of these subspaces is integrable and defines a foliation on $\MS$. In other words:
\begin{dfn}\label{fol} A leaf $\L$ of the foliation on $\MS$ is defined to be the locus along which the periods of the corresponding differentials remain (covariantly) constant.
\end{dfn}
\begin{rem}  In early works \cite{kr-tau,kp1,kp2} (see also \cite{kr-grush1,kr-grush2})
the structure of the foliation was defined on the moduli space $\MO$ of real normalized differentials having poles of exact order $n_\a+1$ at the marked points. The definition above
extends the construction onto the moduli space $\MS$ corresponding to the real normalized differentials with
poles of order at most $n_\a+1$.
\end{rem}
The foliation we define is real-analytic, but its leaves are complex (so the foliation
are real-analytic ``in the transverse direction''). Indeed, locally in the neighborhood $U$ of a curve $\G_0$ one can always choose a basis of cycles on every curve $\G\in U$, which continuously varies with a variation of $\G$. Therefore, locally a leaf $\L$ is defined to be the locus where the integrals of $\Psi$ over the chosen basis of cycles $A_1,\ldots,A_g,B_1,\ldots,B_g$ are equal to $a_1,\ldots,a_g,b_1,\ldots b_g$ --- these are holomorphic conditions, and thus the leaf is a complex submanifold $\L\subset \MS$. If a different basis of $H_1(\G,\mathbb Z)$ is chosen, the periods of $\Psi$ over the basis are still fixed along a leaf (though numerically different).

\begin{rem} In \cite{kr-grush1} the leaves of the foliation on $\MO$ were called big leaves as opposed to small leaves defined by periods of two differentials. They can be regarded as a generalization of the Hurwitz spaces of
$\mathbb P^1$ covers. More precisely, if $k>1$, or if  $k=1$ but $n_1>1$, then there is a special leaf $\L_0\subset \MS$,
on which all the periods of $\Psi$ vanish, i.e., $\Psi$ is exact, $\Psi=dF$. On curves corresponding to points of the open set $\L_0\cap \MO$  the meromorphic function $F$ has poles of orders $n_a$ at $p_\a$, i.e., $F$ defines a cover of $\mathbb P^1$ with prescribed type of branching over one point (infinity). Hence, $\L_0\cap \MO$ can be identified with $C^*$-bundle over the Hurwitz space $\h^{(n_1,\ldots,n_k)}$.
\end{rem}
\begin{theo}\label{coordinates} A leaf $\L$ is a smooth complex subvariety of real codimension $2g$ (\,i.e. of complex dimension $d_g^{\,(n)}=2g-3+|n|+k$).
\end{theo}

\begin{rem}
Theorem \ref{coordinates} extends to the case of $\MS$ the result proved in \cite{kp1} : {\it an open set of $\L$, namely, $\L\cap \MO$ is smooth}.
\end{rem}
\noindent
{\bf Proof.} As it was mentioned above, locally after choosing a basis of cycles on the curves in the neighborhood of $\G_0\in \M_g$, the leaf $\L$ passing through a point $(\G_0,\Psi_0)\in \MS$ is defined by the equations
\beq\label{perab}
\oint_{A_j} \Psi=a_j,\ \ \oint_{B_j}\Psi=b_j,
\eeq
saying that on a curve $\G$ near $\G_0$ there exists a differential with prescribed periods (the reality is implied).
To prove the theorem, one needs to show that these equations are independent.
It turns out that this is indeed the case, and moreover that the set of periods considered as a set of local functions on $\MS$ can be completed to a local coordinate system near any point of $(\G_0,\Psi_0)\in \MS$.

The construction of such local coordinates on $\MO$ is given in \cite{kp1}. A set of holomorphic coordinates on $\L\cap \MO$ is similar to ones in the theory of Hurwitz spaces. They are {\it critical} values of the corresponding abelian integral
\beq\label{E}
F(p)=c+\int^p \Psi,\ \ \ p\in \G,
\eeq
which is a multivalued meromorphic function on $\G$. On $\MO$ the differential has poles of orders $n_\a+1$ at $p_\a$. Therefore, the zero divisor of $\Psi$ is of degree $d_g^{\,(n)}+1$. At the generic point of $\L\cap \MO$, where zeros $q_s$ of $\Psi$ are distinct, the coordinates on $\L$ are the evaluation of $F$ at these critical points:
\beq\label{es}
\varphi_s=F(q_s),\ \ \Psi(q_s)=0,\ \ s=0,\ldots,d_g^{\,(n)},
\eeq
normalized by the condition $\sum_s \varphi_s=0$. Of course, these coordinates depend upon the path of integration needed to define $F$ in the neighborhood of $q_s$. The normalization above is needed to define in addition a common constant $c$ in (\ref{E}). At points of $\L\cap \MO$, where the corresponding differential has multiple zeros $q_{s_1}=\ldots=q_{s_r}$ the local coordinates are symmetric polynomials  $\s_i(\varphi_{s_1},\ldots,\varphi_{s_r}),\, i=1,\ldots,r$ (it is assumed here that the paths of integrations for critical values $\varphi_{s_k}$ are chosen consistently; for more details see \cite{kp1}).

\begin{rem}\label{orderedfs}
A direct corollary of the real normalization is the statement that imaginary parts $f_s={\rm Im}\, \varphi_s$ of the critical values are independent on paths of integration and depend only on labeling of the critical points. They can be arranged into decreasing  order:
\beq\label{order}
f_0\geq f_1\geq \cdots\geq f_{d-1}\geq f_d, \ \ \ d=d_g^{\,(n)}.
\eeq
After that $f_j$ can be seen as a well-defined continuous function on $\MO$, which restricted onto $\L\cap \MO$ is a piecewise  harmonic function. Moreover, as shown in \cite{kr-grush1} the first function $f_0$ restricted onto $\L\cap \MO$ is a {\it subharmonic function}, i.e, it is a function for which the  maximum principle can be applied: {\it $f_0$ has no local maximum on complex subvariety of $\L$ unless it is constant on this subvariety}.
\end{rem}
Now we are going to introduce a set of local coordinates in the neighborhood of a point $(\G_0,\Psi_0)$ of an arbitrary strata $\MOO\subset \MS$. Let $\C$ be the universal curve over $\MS$, i.e. the bundle whose fiber over $(\G,\Psi)\in \MS$ is the curve $\G$ itself. The marked points define sections $s_\a$ of $\C$. A choice of local coordinates $z_\a$ near $p_\a$ on each of the curve $\G$ in the neighborhood of $\G_0$ is equivalent to a choice of the local trivialization of the neighborhood of $s_\a$ in $\C$. Then, the coefficients $r_{\a,i}$ of the expansion of the singular part of $\Psi$ near $p_\a$,
\beq\label{exppsi}
\Psi=\sum_{i=1}^{n_\a}r_{\a,i}z_\a^{-i-1} dz_\a+O(1)dz_\a,
\eeq
can be regarded as local functions on $\MS$.

If $(\G_0,\Psi_0)\in \MOO$, then one can choose the neighborhood $U_{\varepsilon}$ of $(\G_0,\Psi_0)$ in $\MS$, such that within this neighborhood the following inequalities
\beq\label{neigh}
0<c<|r_{\a,m_\a}|,\ \ |r_{\a,i}|<\varepsilon, \ \ i=m_\a+1,\ldots,n_\a\, .
\eeq
hold, where $c$ is a constant.

If $\varepsilon$ is small enough, then the differential $\Psi$ has $d_g^{(m)}+1$ zeros $q_s$ outside the neighborhoods of marked points $p_a$, which tend to the zeros of $\Psi_0$ on $\G_0$,  as $\varepsilon\to 0$.
The remaining $\sum_\a(n_a-m_\a)$ zeros of $\Psi$ tend to the marked points as $\varepsilon\to 0$.
Recall that the critical values of $F$ depend on a choice of the constant $c$ in (\ref{E}).
Under the conventional choice of this constant by the condition $\sum_s \varphi_s=0$ all the critical values might have no limit at $\Mom$. In what follows, the "finite critical values" $\varphi_s,\, s=0,\ldots, d_g^{(m)},$ of $F$ locally in $U_\varepsilon$ to be defined by the normalization $\sum_{s=0}^{d_g^{(m)}}\varphi_s=0$.

Let us introduce the following set of functions $\{x_A\}$ in $U_\varepsilon \subset \MS$. They are:

$(i)$ the leading coefficients $r_{\a,i},\, i=m_\a+1,\ldots,n_\a,$ of the expansion (\ref{exppsi});

$(ii)$ the finite critical values $\varphi_s, \, s=1,\ldots, d_g^{\,(m)}$ of $F$, if $\Psi_0$ has simple zeros. If $\Psi_0$ has a multiple zero $q_{s_1}=\ldots,q_{s_r}$, then the corresponding subset of finite critical values of $F$ should be replaced by their symmetric polynomials.

\begin{lem}\label{coord1} Near each point $(\G_0,\Psi_0)\in \MO\subset \MS$ the periods $a_j,b_j,j=1,\ldots g,$
of \,$\Psi$, (\ref{perab}), and real and imaginary parts of the functions $x_A$ ((i,ii) above) have linear independent differentials and thus define a real analytic local coordinate system on $\MS$.
\end{lem}
The proof of Lemma \ref{coord1} is similar to that of Theorem 1 in \cite{kp1}.
We outline the key step in it. Suppose that the differentials of these functions are linearly dependent at $(\G_0,\Psi_0)$ (and thus the functions do not give local coordinates near
$(\G_0,\Psi_0)$). Then there exists a one-parametric family $(\G_t,\Psi_t)\in\MS$, with real
parameter $t$ such that the derivative of any of the above functions along this family
is equal to zero at $t=0$.

Recall, that locally on each of the curve $\G_t$ we already fixed a basis $A_j,B_j$ for cycles
(needed for the definition of periods (\ref{perab})). Let $\omega_j(t)$ be the basis
of holomorphic differentials on $\G_t$ dual to $A_j$. Denote then
$v_j(p,t):=\int_{q_1(t)}^p\omega_i(t)$ the corresponding abelian integral --- the
function of $p\in\G_t$, depending on the choice of the path of integration. Let us
also denote $F_t(p):=\int_{q_1(t)}^p\Psi_t$ the integral of our chosen meromorphic
differential along the same path.

We will now want to see how $v_j$ varies in $t$. For this to make sense as a partial
derivative, we need to ``fix'' how to vary the point $p$ as we vary $\G_t$, and to do this we will
use $F_t$ as the local coordinate on the universal cover of $\G$. This is to say that we
will fix $F=F_t(p_t)$ and let $t$ vary; this is to say that $F$ allows us to define a
connection on the space of abelian integrals.

Rigorously, we consider the {\it partial derivative with fixed $F$}:
\begin{equation}\label{partial}
\p_t v_j(F):=\frac{\p}{\p t}v_j(F_t^{-1}(F),t)|_{t=0}
\end{equation}
and show that it is zero.
We think of the surface $\G_t$ as cut along a basis of
cycles, so that the integrals $\int_{q_1(t)}^p$ in the definition above are taken
along paths not intersecting this basis, i.e.~on the simply-connected cut surface.
Then the expression above is by definition holomorphic on the cut surface $\G_0$, with
possible poles at the zeros of $\Psi$ (where $F^{-1}$ is singular), and with
discontinuities along the cuts. However, if as we wary $t$
the periods of $\Psi_t$ do not change, (\ref{partial}) has constant "jumps" along the cuts, and since the critical values $\varphi_s$ do not change, (\ref{partial}) also has no poles at zeros  of $\Psi$. By chain rule the partial derivative with fixed $F$ and partial derivative with {\it fixed $z_a$} (which provides a trivialization of the neighborhood of
the marked point $p_{\a,t }\in \G_t$) are related by the equation
\beq\label{partialz}
\p_t v_j(F)=\p_t v_j(z_\a)-\frac{\omega_j}{\Psi_0}\,\p_t F(z_\a)|_{\, t=0}.
\eeq
By our assumption the first $(n_\a-m_\a)$ leading coefficients of the expansion (\ref{exppsi}) do not change.  Then from (\ref{partialz}) it follows that (\ref{partial}) has no pole at $p_\a$.
Thus the differential of expression (\ref{partial}) is a holomorphic differential on
$\Gamma_0$ with zero $A$-periods; thus it is identically zero, and also has zero
$B$-periods. Since the  $B$-periods of $\omega_j$ are entries of the period matrix
$\tau$ of $\G_0$, this means that we have
$$
 \frac{\p}{\p t}\tau_{ij} (t)|_{t=0}=0,
$$
The infinitesimal Torelli theorem says that the period map $\tau:\M_g\to\A_g$ induces
an embedding on the tangent space away from the locus of the hyperelliptic curves and
thus the above is impossible unless $\G_0$ is a hyperelliptic curve.
For a hyperelliptic curve the kernel of $d\tau$ is one-dimensional and transverse to
the tangent space of the hyperelliptic locus.
Therefore, to complete the proof, it suffices to show that if $\G_0$ is a hyperelliptic curve, then the tangent vector to the family $\G_t$ at $t=0$ is tangent to the locus of hyperelliptic curves.
The proof of this last step is identical to that in \cite{kp1}.

\begin{cor} The set of functions $x_A$ define a system of local complex coordinates
on the leaf $\L$ passing through the point $(\G_0,\Psi_0)$.
\end{cor}

\medskip
For further use,  we consider now  a partial compactification of $\MS$. The space of real normalized differential is invariant under the multiplication by {\it real numbers}. Let $\PS=\MS/R_+$ be the factor-space under the action of the multiplicative group of positive real numbers.
The fiber of the forgetful map $\PS\longmapsto \M_g$ can be seen as a space of "unit" singular parts. It is isomorphic to the sphere $S^{2|n|-1}$.

The canonical foliation of $\MS$ descents to a foliation structure on $\PS$.
\begin{dfn} A leaf $[\L]$ of the foliation on $\PS$ is defined to be the projection of the locus in $\MS$ along which the ratio of any two periods of the corresponding differentials remains constant (if two periods are zero, they both must remain zero).
\end{dfn}
The multiplication by real numbers acts "transversally" on all leaves $\L$ of the big foliation except the Hurwitz leaf $\L_0$ corresponding to exact differentials. Therefore, a leaf
$[\L]$ passing through a point $(\G,[\Psi])\in \PS$, which is not in the image of $\L_0$, is locally isomorphic  to the leaf $\L\in \MS$ passing through a point $(\G,\Psi)\in \MS$ in the preimage of
$(\G,[\Psi])$. Hence, $[\L]$ has a natural complex structure. (Notice that we can not say anymore that $[\L]$ is a local complex subvariety because $\PS$ has no complex structure by itself.)

The leaf $[\L_0]$ is the only singular leaf of the foliation. It is not complex, and has real dimension one less that all the other leaves. It is isomorphic to a $S^1$-bundle over the Hurwitz space.

\section{Dual cycles and periods}

In this section we introduce an additional tool needed for the proof of Arbarello's conjecture in the next section: a notion of cycles dual to critical points of real normalized differentials.

For simplicity, from now on we consider only the case $(\G,\Psi)\in \Mo$ of real normalized differentials having pole only at one marked point ($k=1$).
By definition of the real normalization, the imaginary part $\Phi={\rm Im}\, F$ of an abelian integral $F$
of such differential is a single-valued harmonic function on $\G\setminus p_1$. The level curves  $\Phi_h:=\{p\in \G: \Phi(p)=h\}$ of this function is a union of cycles that are smooth everywhere except at $p_1$ and at $q_s$ if $h=f_s=\Phi(q_s)$. For $h$ big enough ($h>f_0)$ the level curve is a union of $n$ "loops" in the small neighborhood of $p_1$. The real part of $F(p)$ is multi-valued. Nevertheless, everywhere except, at zeros of $\Psi$, the direction along which the real part remains (locally) constant is well-defined. Let us call the integral lines of these directions {\it imaginary} rays. It will be always assumed that they are oriented such that $\Phi$ increases along the orientation. In the small neighborhood of $p_1$ one can always fix a single-valued branch of $F$. If the imaginary ray passing through a point $p$ does not pass through the zeros of $\Psi$, then the real part of $F$ evaluated at points of the ray within the neighborhood of $p_1$ uniquely defines the real part of $F(p)$. That allows one to define $F(p)$ as a single-valued holomorphic function on $\G\setminus \Sigma$, where $\Sigma$ is a graph, whose edges are imaginary rays that begin at $p_1$ and zeros of $\Psi$, and end at zeros of $\Psi$. By continuity $F$ can be extended on $\G$ cut along the edges of $\Sigma$. But the limits $F^{\pm}(p), \, p\in \Sigma$ on two sides of the cuts are in general distinct. The discontinuity of $F$, i.e. the "jump" function $j(p):=F^+(p)-F^-(p)$ is constant on each of the edges.

Let $\D$\, be an open subset of $\Mo$, where the differentials have simple zeros with "distinct" real parts of the critical values. The later means that imaginary rays emanated from a zero of the differential do not contain any other zero, i.e. they end at $p_1$ . In this case the graph has $2g+n-1$ connected components $\Sigma_s$. Each component is a union of two imaginary rays that both begin at $p_1$ and end at one of the zeros $q_s$ of $\Psi$. Along these rays the function $\Phi$ increases from $-\infty$ through $f_s$, i.e. each of the zeros $q_s$ of $\Psi$ is a "tip" of $\Sigma_s$. At the same time, for each zero $q_s$ there are two imaginary rays that both begin at $q_s$ and end at $p_1$. Along these rays $\Phi$ increases from $f_s$ through $+\infty$.
Reversing the orientation of one of these rays we can define a closed oriented cycle $\sigma_s$.

Although the differential $\Psi$ has singularity on $\s_s$, the period of $\Psi$ over the homology class $[\s_s]\in H_1(\G,\mathbb Z)$ of $\s_s$ is well-defined (recall that $\Psi$ has no residue at $p_1$). It equals
\beq\label{periodsigma}
\pi_s:=\oint_{[\,\s_s]}\Psi=r_s^1-r_s^2,\
\eeq
where $r_s^{1}, r_s^2$ are two values of the real part of $F(q_s)$ asymptotically defined along each of the rays. Notice, that on the edges of the graph $\Sigma$ the jump function equals $j(p)=\pm \pi_s,\ \ p\in \Sigma_s$.

\begin{rem}\label{n1} In the simplest case of real normalized differentials having one pole of the second order ($n=1$),
the construction of dual cycles  looks especially tempting. In that case, if $\Psi$ has simple zeros with
distinct real parts of the critical values, then it defines on the corresponding curve $2g$ dual cycles. It is easy to see that each of the dual cycles has non-trivial homology class. Indeed, for $n=1$ the abelian integral $F$ of $\Psi$ has simple pole at $p_1$. That implies that the period of $\Psi$ over $[\s_s]$ is never zero, $\pi_s\neq 0$. Hence, $[\s_s]\neq 0$, as well. Moreover, it turns our $[\s_s]$ are linear independent and, thus, define a basis in $H_1(\G,\mathbb Z)$ (see Lemma \ref{span} below).
\end{rem}
At points $(\G,\Psi)$, where $\Psi$ has multiple zero, or where there are critical values of $F$ with
"non-distinct" real parts, the structure of the graph $\Sigma$ might be combinatorially relatively complex. Still, for each zero $q_s$ of $\Psi$ there is only a {\it finite} number  of semi-infinite paths $\g_s^i$  (arbitrary ordered)
along imaginary directions which start at $q_s$ and end at $p_1$.  A pair of such paths define an oriented cycle $\s_s^{ij}=\g_s^i\cup (-\g_s^j),\, i<j$.
\begin{lem}\label{span} The homology classes $[\s_s^{ij}]$ of dual cycles  span the homology group $H_1(\G,\mathbb Z)$.
\end{lem}
The proof of the lemma is just another incarnation of the arguments used in
\cite{kp1} and slightly extended in the proof of Theorem \ref{coordinates}.
Indeed, suppose that the classes $[\s_s^{ij}]$ defined by the differential $\Psi_0$ on a curve $\G_0$ do not span $H_1(\G_0,\mathbb Z)$, then by
Theorem 1 in \cite{kp1}, there is a one-parametric deformation $(\G_t,\Psi_t)\in \Mo$ such that
that the derivatives of critical values $\varphi_s$ and the periods $\pi_s^{ij}$ along this family are equal to zero at $t=0$.

Recall once gain, that locally it can always be assumed that on each of the curve $\G_t$
one has a fixed basis $A_j,B_j$ for cycles. Let $\omega_j(t)$ be the basis
of holomorphic differentials on $\G_t$ dual to $A_j$. Let us fix a branch of the corresponding abelian integral $v_j$ in the sectors of the neighborhood $p_1$ where $\Phi\to +\infty$. Then extending it along imaginary rays (backward) we can define a single-valued branch of $v_j$ on
$\G_t\setminus \Sigma_t$, and then extend it by continuity on two sides of edges of $\Sigma_t$. Jumps of $v_j$ on edges of the graph are linear combinations of periods of $\omega_j$ over dual cycles.
As before, we then consider the partial derivative $\p_t v_j(F,t)$ with fixed $F$ at $t=0$. If the derivatives of $\varphi_s$ and jumps of $F$ on the edges of $\Sigma$ are equal to zero (as assumed), then
the derivative $\p_t v_j(F,t)|_{t=0}$ is holomorphic on $\G_0\setminus \Sigma_0$, and has constant jumps on $\Sigma_0$. Then, the differential of that expression is a holomorphic differential on $\G_0$ with zero $A$-periods. Hence, it is identically zero. That implies that the matrix of $b$-periods of $\omega_j$ do not change  along the family at $t=0$. The remaining steps of the proof are identical to that in \cite{kp1}.

Let $r_s^i$ be "real values" of $F(q_s)$ defined asymptotically by the paths $\g_s^i$ .
In the general case, the jumps $j(p)$ on edges of $\Sigma$ are linear combinations with integer coefficients
of periods of $\Psi$ over the dual cycles
\beq\label{periodsigma1}
\pi_s^{ij}:=\oint_{[\,\s_s^{ij}]}\Psi=r_s^i-r_s^j.
\eeq

\section{Proof of Arbarello's conjecture.}

For the motivation of the further steps towards the proof of Arbarello's conjecture it is instructive to outline the new proof of Diaz'  theorem via real-normalized differentials with one pole of the second order.

\medskip
\noindent{\bf Diaz'\,  theorem revisited.} Let $\phi_s$ be a "weighted critical value" of $F$ defined by the formula
\beq\label{weightphi}
\phi_s=\frac{\varphi_s}{|\,\pi_s|}\,,\ \ |\,\pi_s|:=\min_{ij}\, \{|\,\pi_s^{ij}|\neq 0\},
\eeq
where the minimum is taken over a finite set of nonzero periods dual to the critical value.
As it was stressed in Remark \ref{n1}, dual periods of a real normalized differential with one pole of second order never vanish. Therefore, $\phi_s$ is a well-defined local function on $\Moo$.
Imaginary parts of the weighted critical values $g_s={\rm Im}\, \phi_s$ are independent of paths of integration and depend only on the labeling of the critical point. They can be arranged
into decreasing order:
\beq\label{orderedgs}
g_0\geq g_1\geq \cdots\geq g_{2g-2}\geq g_{2g-1},
\eeq
After that $g_s$ can be seen as a well-defined function
on $\Moo$ (compare with Remark \ref{orderedfs}). It is continuous on the open set $\D\subset \Moo$, where the corresponding differentials $\Psi$ has simple zeros with distinct real part of the critical values.
(Recall, the later means that imaginary rays emanated from a zero of the differential do not contain any other zero.)
Moreover, $g_s$ restricted onto $\L\cap \D$ is a piecewise harmonic function.

From (\ref{weightphi}) it is easy to see that the evaluation of $g_s$ at any point of $\Moo$
is equal to the maximum of the limits at this point of values of $g_s$ on $\D$. Therefore,
$g_s$ is {\it upper semicontinuous} function on $\Moo$.
\begin{lem}\label{g0} The function $g_0$ restricted onto a complex subvariety of $\L$ has a local maximum if and only if it is constant along the subvariety.
\end{lem}
{\bf Proof.} By definition of $\L$, the periods $\pi_s$ are constant on $\L\cap \D$, therefore on
$\L\cap \D$ the function $g_0$, as the maximum of harmonic functions, has local maximum only if it is a constant. In order to show that $g_0$ have no  local maximum at a point of discontinuity,
it is enough to notice, that the "direction" of its discontinuity, which correspond to a change
of the real part ${\rm Re}\, \varphi_0$ is always transversal to the directions along which the imaginary
part ${\rm Im}\, \varphi_0$ changes.

\medskip
Let $X$ be a compact cycle in $\M_g$. Its preimage $Y$ under the forgetful map $\Moo\longmapsto \M_g$ is not-compact, but the factor-space $Z=Y/R_+$ is compact. The function $g_0$ is homogeneous under multiplication of real normalized differentials by a positive real number, i.e. it descents to a upper semicontinuous function $g_0'$ on $\Poo$. Therefore, there is a point $(\G_0,[\Psi_0])$ of $Z$, where the restriction of $g_0'$ on $Z$ achieves its supremum.
Let us fix a preimage $(\G_0,\Psi_0)\in \Moo$ of this point. At this point the function $g_0$ achieves its supremum on $Y$.

Let $\L$ be the leaf of the big foliation passing through $(\G_0,\Psi_0)\in \Moo$.
At this point the function $g_0$ restricted onto $\L\cap Y$ has a local maximum. Then, by Lemma \ref{g0} it must be a constant on $\L\cap Y$.

Let $Y_0\in Y$ be a preimage of the compact set $Z_{0}\subset Z$, where $g'_{0}$ takes its maximum value. On the compact set the second function $g_1'$ must achieve its supremum.
As it was shown above, $Y_0$ is foliated by leaves $\L\cap Y$. On these leaves the second function $g_1$ is subharmonic, i.e., it must be a constant. Continuing by induction we get that all the functions $g_s$ are constants on $\L\cap Y$. If $g_s$ are constants, then $\varphi_s$ are (locally) constants on $\L\cap Y$, as well. The functions $\varphi_s$ are local coordinates on $\L$. Therefore, $\L\cap Y$ must be zero-dimensional. But, if $X$ is of dimension greater that $g-2$, then $\L\cap Y$ is at least one-dimensional. The contradiction we get completes the proof of Diaz' theorem.

\medskip
The proof above is almost a carbon copy of the proof of Diaz' theorem in \cite{kr-grush1} with the only (but important) modification: the replacement of critical values by weighted critical values. The proof of Arbarello's conjecture follows mainly the same line of arguments but requires further modifications due to the following: (a) periods over dual cycles of a real normalized differential having pole of order $n>1$ might vanish, (b) only part of critical values of such differential have finite limits on smaller strata $\Mom$. Notice, that these two issues
are interconnected, at least in the neighborhood $U_{\varepsilon}$ of $\Mom,\, m<n$. Indeed, for sufficiently small $\varepsilon$ periods over the cycles dual to the critical points that are in a small neighborhood of $p_1$ are always zero, because the cycles are homologically trivial.

For further use, let us introduce on $\Mon$ related strata. First, for $(\G,\Psi)\in \Mon$ consider a subset $S$ of critical points of $\Psi$ such that at least one of the periods of $\Psi$ over the dual cycles $\s_s^{ij}$ is not zero, i.e. $\{q_s\in S: \ \exists\ \,  \pi_s^{ij}\neq 0\}$.

\begin{dfn} The stratum $\SP_k$ of $\Mon$ is defined to be the locus, where $\Psi$ has exactly $k$ zeros with non-vanishing periods over the dual cycles, $|\,S|=k$.
\end{dfn}

At points of $\SP_k$ let us normalize the critical values of the corresponding differential by the condition $\sum_{q_s\in S}\varphi_s=0$. Then, as before, the imaginary parts ${\rm Im}\, \phi_s, \, s\in S,$ of the weighted critical values (\ref{weightphi}) can be arranged into decreasing order:
\beq\label{wns}
w_0^k\geq \ldots\geq w_{k-1}^k.
\eeq
After that $w_j^k$ can be seen as well-defined upper semicontinuous function on $\SP_k$.
The functions $w_j^k$ are continuous on $\SP_k\cap \D$ and $w_0^k$ is subharmonic on $\SP_k\cap \L$.

On a leaf $\L$ periods of $\Psi$ over any continuously varying cycles are constants. Therefore, the locally constant function $|S|$ is discontinuous only at points, where one of the zeros $q_s$ of the corresponding differential is excluded from $S$ under variation moving it out of some edge of the graph $\Sigma$. Notice, that if it is not in $S$ on one side of the edge, then it remains in $S$ under variations moving it on another side of the edge. Therefore,
the intersection $\L\cap \overline{\SP_k}$ of a leaf with the closure of $\SP_k$ is a complex domain with boundaries. Notice also, that the imaginary part of the vanishing critical value should be less than the imaginary part of the critical value corresponding to the top vertex of the graph edge. At the same time, the set of not trivial periods dual to the vanishing critical value is a subset of the periods dual to the top vertex of the edge. Hence, vanishing weighted critical value satisfies the inequality ${\rm Im}\, \phi_s\leq w_0^k$. Thus,  at common points of the boundaries of $\L\cap \overline{\SP_k}$ and
$\L\cap \overline{\SP_{k+1}}$ the functions $w_0^k$ and $w_0^{k+1}$ coincide, i.e., on each leaf $\L$ of the foliation
there is a well-defined subharmonic function $W_0^{\L}$. It is necessary to emphasize that the arguments above hold only after restriction of consideration onto the leaves of the foliation. Under variation transversal to the leaves the vanishing critical might tend to infinity, and therefore, there is no globally defined upper semicontinous function $W_0$.

Now we are ready to present the proof of the main theorem.

\begin{theo} Any compact complex cycle in $\M_g$ of dimension $g-n$ must intersect the locus $\W_n\subset \M_g$
of smooth genus $g$ algebraic curves having a Weierstrass point of order at most $n$.
\end{theo}
\noindent{\bf Proof.} Let $X$ be a complex compact cycle in $\M_g$ that does not intersect $\W_n$. Then its preimage $Y$ under forgetful map $\Mon\longmapsto \M_g$ does not intersect the leaf $\L_0\subset \Mon$ --- the locus of exact differentials. Therefore, for each point $(\G,\Psi)\in Y$ the corresponding set $S$ of critical points is always non-empty. Let $k_0$ be the minimum of $|S|$ on $Y$, i.e. the minimum number such that $Y^k=\SP_k\cap Y$ is not empty.

In general the loci $Y^k$ are not closed, partly, because of vanishing of some periods at limiting points. Although $k_0$ is the minimal number of critical values with non-zero dual periods, the locus $Y^{k_0}$ still might be not closed. Points of the closure $\overline {Y^{k_0}}$ that are not in $Y^{k_0}$ are points in $Y^k$ with $k>k_0$, where one or more critical points vanish under variations along leaves of foliation. At these points the limiting values of the functions $w_j^{k_0}$ are finite and are a subset of functions $w_j^k$. Hence, the functions $w_j^{k_0}$ defined first
on $Y^{k_0}$ extend, as upper semicontinuous functions, on the closure $\overline{Y^{k_0}}$.

The factor-space $\overline{Z^{k_0}}=\overline{Y^{k_0}}/R_+$ is compact. The function $w_0$ descents to an upper-semicontinuous function $w_0'$ on $\overline{Z^{k_0}}$, which must achieve its supremum at some point. Let $\L$ be the leaf of the foliation passing through any preimage of that point on $\overline{Y^{k_0}}$. Then $w_0$ onto $\L\cap \overline{Y^{k_0}}$ has a local maximum. The intersection $\L\cap \overline{Y^{k_0}}$ is a complex domain with boundary. The arguments used in the proof of Lemma \ref{g0} show that $w_0^{k_0}$ is subharmonic in the interior points of the domain. The similar arguments show that any value of $w_0^{k_0}$ at a point of the boundary is achieved at some interior point. Thus $w_0^{k_0}$ must be a constant. If it is a constant, then the next function $w_1^{k_0}$ on $\L\cap Y^{k_0}$ is subharmonic and by the same arguments must be a constant. Continuing by induction we get that all the functions $w_s$ are constants. If $w_s^{k_0}$ are constants,  then all the critical values $\varphi_s, \, q_s\in S,$ are constants, as well.

Consider now the locus $Y_0\subset Y$, where in the set of all the critical values there exists $k_0$-tuple of zeros of $\Psi$ such that: (a) the ordered tuple $w_0\geq\ldots\geq w_{k_0}$ of the imaginary parts of the weighted critical values remain constant along the corresponding leaf of the foliation, and (b) the imaginary part of the first critical value is normalized by the condition $f_0=1$. Notice, that the later normalization fixes the lifting of $Y_0/R_+\subset Z=Y/R_+$ into $Y$.

The locus $Y_0$ is compact (and as shown above is non empty). Therefore, the continuous function
$\wt f={\rm max}_{s\in S} \,{\rm Im}\, \phi_s$ restricted onto $Y_0$ achieves its supremum. (Recall, that by definition the zero $q_s$ of the differential $\Psi$ with $s\in S$ can not approach $p_1$, and therefore,
the function $\wt f$ is bounded on $Y$. The function $\wt f$  restricted onto any leaf of the foliation is subharmonic. Hence, using the same arguments as above, we get that on $Y_0\cap \L$
all the critical values having non-zero dual period are constant.

Our next goal is to show that all the remaining critical values $\varphi_s, q_s\notin S,$ are also constants on $\L\cap Y_0$.  Suppose not, then the "moving" along $\L\cap Y_0$ critical point $q_s\notin S$ must intersect
a "non-moving" edge of the graph $\Sigma$ with a non-trivial jump on it. After such crossing at least one period  dual to $q_s$ becomes non-zero. If $\phi_s$ becomes  a constant after crossing the edge of the graph, then it should be a constant before the crossing. Thus, we get that all the critical values are constant along $\L\cap Y_0$.
By definition of $Y_0$, we have $\L\cap Y_0=\L\cap Y$. If $\varphi_s$ are constants, then $\L\cap Y$ is at most zero-dimensional. That contradicts to the assumption that $X$ has dimension at least $g-n$. The theorem is proved.


\begin{thebibliography}{}
\bibitem{arb}
E.~Arbarello, {\it Weierstrass points and moduli of curves}, Compositio Math. {\bf 29} (1974), 325--342.
\bibitem{arbmond} E.~Arbarello, G.~Mondello, {\it Two remarks on the Weierstarss flag}, in: "Compact Moduli spaces and Vector Bundles", ed.
Valery Alexeev, Angela Gibney, Elham Izadi and David Swinarski, Contemporary Mathematics (Proceedings) series, vol.564. (to appear)
\bibitem{brad} H.W.~Braden, I.M.~Krichever (eds.), {\it Integrability: The Seiberg-Witten and Whitham equations}, Amsterdam: Gordon and Breach Science Publishers. 2000,
\bibitem{diaz}S.~Diaz, {\it Exceptional Weierstrass points and the divisor on moduli space that they define.} Mem. Amer. Math. Soc., {\bf 56} no. 27. Providence, RI, 1985.
\bibitem{gmmmk} A.~Gorsky, I.~Krichever, A.~Marshakov, A.~Mironov, A.~Morozov, {\it Integrability and Seiberg-Witten Exact Solution}, Phys.Lett. {\bf 355B} (1995), 466-474
\bibitem{har-mor} J.~Harris, I.~Morrison, {\it Moduli of curves}, Graduate Texts in Mathematics, 187. Springer-Verlag, New York, 1998.
\bibitem{kr-grush1}
S.~Grushevsky, I.~Krichever, {\it The universal Whitham hierarchy and geometry of the moduli space of pointed Riemann surfaces}, in Surveys in Differential Geometry {\bf 14} (2010).

\bibitem{kr-grush2} S.~Grushevsky, I.~Krichever, {\it  Foliations on the moduli space of curves, vanishing in cohomology, and Calogero-Moser curves}, arXiv:1108.4211

\bibitem{kr-real} I.~Krichever, {\it The spectral theory of ``finite-gap'' nonstationary Schr\"odinger operators}, The nonstationary Peierls model. (Russian) Funct. Anal. Appl. 20 (1986), no. 3, 42--54, 96.

\bibitem{kr-av} I.~Krichever: {\it Averaging method for two-dimensional integrable equations}, Funct. Anal. Appl. {\bf 22} (1988), 37--52.

\bibitem{kr-tau} I.~Krichever: {\it The $\tau$-function of the universal Whitham hierarchy, matrix models, and topological field theories}, Comm. Pure Appl. Math. {\bf 47} (1994), 437--475.

\bibitem{kp1} I.~Krichever, D.H.~Phong: {\it On the integrable geometry of $N=2$ supersymmetric gauge theories and soliton equations}, J. Differential Geometry {\bf 45} (1997) 445--485.
\bibitem{kp2} I.~Krichever, D.H.~Phong: {\it Symplectic forms in the theory of solitons}, Surveys in Differential Geometry {\bf IV} (1998), edited by C.L. Terng and K. Uhlenbeck, 239-313, International Press.
\bibitem{lo}E.~Looijenga: {\it On the tautological ring of ${\M}_g$}, Invent. Math.  {\bf 121}  (1995), no. 2, 411--419.


\end{thebibliography}
\end{document}